\documentclass[12pt]{article}
\usepackage{amssymb}
\usepackage{epic}
\usepackage{eepic}
\usepackage[dvips]{graphicx}
\usepackage{epsfig}

\begin{document}

\def\stackunder#1#2{\mathrel{\mathop{#2}\limits_{#1}}}
\def\binom#1#2{{#1 \choose #2}}

\newtheorem{theorem}{Theorem}[section]
\newtheorem{proposition}[theorem]{Proposition}
\newtheorem{definition}[theorem]{Definition}
\newtheorem{corollary}[theorem]{Corollary}
\newtheorem{lemma}[theorem]{Lemma}
\newtheorem{remark}[theorem]{Remark}
\newtheorem{example}[theorem]{Example}

\title{{\Huge{\textbf{Knots, Feynman Diagrams and
Matrix Models}}}}
\author{\textbf{Martin Grothaus} \\
BiBoS, Universit\"at Bielefeld, D 33615 Bielefeld,
Germany \and
\textbf{Ludwig Streit} \\
BiBoS, Universit\"at Bielefeld, D 33615 Bielefeld,
Germany \\
CCM, Universidade da Madeira, P 9000 Funchal, Portugal
\and
\textbf{Igor V.~Volovich} \\
Steklov Math.~Inst., RAS, R 117966 Moscow, Russia}
\date{\today}
\maketitle

\bigskip

\begin{abstract}
An $U(N)$-invariant  matrix model with $d$ matrix variables
is studied. It was shown that in the
limit $N \to
\infty$ and $d \to 0$ the model describes the knot
diagrams.
We realize the free partition function of the matrix model as the generalized
expectation of a Hida distribution $\Phi_{N, d}$. This enables us to give a
mathematically rigorous meaning to the partition function with
interaction. For the generalized function
$\Phi_{N, d}$ we prove a Wick theorem and we derive explicit formulas
for the propagators. 
\end{abstract}

\pagebreak

\tableofcontents

\pagebreak

\section{Introduction}

In this note we study an $U(N)$-invariant matrix model with the
following partition function
\begin{eqnarray}\label{eq003}
Z(N, d, g) = \int_{{\mathbb R}^{2dN^2}} \exp(iS(A,B,g)) \,dA \,dB,
\end{eqnarray}
where the action is
\begin{eqnarray*}
S(A, B, g) = Tr(A_{\mu}B_{\mu}) + \frac{g}{2N}
Tr(A_{\mu}B_{\nu}A_{\mu}B_{\nu}). 
\end{eqnarray*}
Here $N, d \in {\mathbb N}$, $A_{\mu}$ and $B_{\mu}$ are Hermitian $(N
\times N)$-matrices over the field ${\mathbb C}$, $1 \le \mu \le d$. 
We assume the summation over the
repeating indices and $g \in {\mathbb R}$ is the coupling constant. The measure
\begin{eqnarray*}
dA := \prod_{\mu = 1}^{d} \Bigg{(} (\prod_{1 \le k \le l \le N} \,d
\Re{A}^{kl}_{\mu}) (\prod_{1 \le k < l \le N} \,d\Im{A}^{kl}_{\mu}) \Bigg{)}.
\end{eqnarray*}
Our investigations are motivated by a paper of Aref'eva and Volovich \cite{AV98} in which the
authors proved the following theorem.
\begin{theorem}\label{th1000}
The set of connected Feynman diagrams for the model (\ref{eq003}) in the
limit $N \to \infty$ and $d \to 0$ is in one-to-one correspondence
with the set of alternating knot diagrams. The generating function for
the alternating knot diagrams is given by the expression
\begin{eqnarray*}
F(g) = \lim_{d \to 0} \lim_{N \to \infty} \frac{1}{dN^2} \ln Z(N, d ,g).
\end{eqnarray*}
\end{theorem}

Theorem \ref{th1000} can be regarded is an application of quantum
field theory to the theory of knots. 
Another approach has been performed by Witten
\cite{Wi89}. There the author has considered 
the Wilson loops in Chern-Simons gauge theory as knots. 
In general the path integrals used in quantum field theories are not
mathematically well-defined. In 
\cite{AS92} and \cite{Sc91} a rigorous mathematical
meaning has been given to the path integral of the Abelian
Chern-Simons 
model, using the theory of Fresnel integrals developed in \cite{AlHK76};
and in \cite{LS95} this has been done in the frame of white noise
analysis, see \cite{HKPS93}. Furthermore, in
\cite{LS95} the authors have specified the Wilson loops and have proved
the relation to topological invariants which has been conjectured in
\cite{Wi89}. In the non-Abelian case a mathematically rigorous meaning to
the Chern-Simons integral has been 
given in \cite{AS96}, also by using the white noise approach.

The theory of knots, see \cite{At90} and \cite{Ka93}, is used in low
dimensional topology and also in physics, chemistry and,
biology. Recently Faddeev and Niemi \cite{FN97} have suggested that in
certain relativistic quantum field theories knot-like configurations
may appear as stable solitons. The remarkable progress in the
classification of knots, based on Jones' and Vassiliev's
invariants, has been related with the application of 
von Neumann algebras, Yang-Baxter equations, and singularity theory,
for a review see \cite{Bi93} and \cite{BGRT97}. 

The large $N$ limit is considered in quantum chromodynamics (QCD), matrix
models and superstring theory, see \cite{Ho74}, \cite{BIPZ78},
\cite{AAV95}, \cite{AV96}, 
\cite{BFSS96}, \cite{IKKT96}, and \cite{DVV97}, the limit $d \to 0$ is
considered in the theory of spin glasses and in polymer physics, see
\cite{Ge82} and \cite{Ar83}. 
 
The statement of Theorem \ref{th1000} reminds the matrix approach to
superstring theory, see \cite{BFSS96}, \cite{Pe96}, \cite{IKKT96}, and
\cite{DVV97}, where space-time is represented as the moduli space of
vacuum and strings appear in the large $N$ limit. See also the
concluding remarks given in Section \ref{ss4711}.

In order to prove Theorem \ref{th1000}, in \cite{AV98} the authors have expanded
the partition function $Z$, see (\ref{eq003}), into the formal
perturbation series in the coupling constant $g$ 
\begin{eqnarray}\label{eq004}
Z = \sum_{k = 0}^{\infty} \frac{1}{k!} \Big{(} \frac{ig}{n}
\Big{)}^k \int_{{\mathbb R}^{2dn^2}} (Tr(A_{\mu} B_{\nu} A_{\mu}
B_{\nu}))^k \exp(iTr(A_{\mu} B_{\mu})) \,dA \,dB. 
\end{eqnarray}
Our contribution to the proof of Theorem \ref{th1000} is to give a
mathematically rigorous meaning to the Gaussian integrals in (\ref{eq004}), to
prove a Wick theorem for these integrals, and to derive explicit
formulas for the propagators. 

This paper is organized as follows. In Section
\ref{s123} we introduce
the concepts of Gaussian analysis as far as necessary
for the
constructions in this note. In doing so we concentrate
on the special case
in which the underlying nuclear triple is given by
matrix spaces. In
order to give a mathematical rigorous meaning to the
integrals in
(\ref{eq004}) the spaces of generalized functions
(Hida distributions)
we introduce in Section \ref{ss23} are essential. For
a detailed exposition of
Gaussian and white noise analysis we refer to the
monographs
\cite{Hi80}, \cite{BeKo88}, \cite{HKPS93},
\cite{Ob94},
and \cite{Kuo96}.
In Section
\ref{s0815} we realize the partition function (\ref{eq003}) for $g = 0$ as the
generalized expectation of a Hida distribution $\Phi_{N, d}$. At first we consider
the one dimensional case, i.e., $N = 1$, and then the case of an
arbitrary $N \in {\mathbb N}$. In our first approach we obtain for different $N$
different spaces of Hida distributions. In order to calculate the large $N$
limit it turns out to useful to embed all $\Phi_{N, d}$ in a joint
space of Hida distributions. This embedding we derive in Section
\ref{ss626}. In Section \ref{s777} we discuss the proof of Theorem
\ref{th1000}. First we shortly recall the necessary definitions and
notions from the theory of knots. In Section \ref{ss1860} we give a
mathematically rigorous definition of the Gaussian integrals in
(\ref{eq004}), see Definition \ref{de1001}, prove a Wick theorem for the
generalized function $\Phi_{N, d}$, see Theorem \ref{th1003}, and
derive explicit formulas for the propagators, see Theorem
\ref{pr1000}. Then we show how to apply the 
Feynman diagram technique in order to find an explicit formula for the
logarithm of the partition function, see (\ref{eq357}). This formula enables us to
determine the limit $N \to \infty$ and $d \to 0$ of the normalized
logarithm of the partition function and to prove Theorem
\ref{th1000}. Finally, in Section \ref{ss4711} we quote some
concluding remarks on the consequences of Theorem \ref{th1000} for
M(atrix)-theory given in \cite{AV98}.  

\section{Gaussian analysis on matrix
spaces}\label{s123}

\subsection{Gaussian spaces on matrix spaces}

We start by considering a Gel'fand triple
\begin{eqnarray*}
{\cal M} \subset {\cal T} \subset {\cal M}^\prime.
\end{eqnarray*}
Here $\cal T$ is a real separable Hilbert space of
matrices. Its inner product $(\cdot,\cdot)$ we assume
to be given by
the trace of the matrix product of two matrices and
the corresponding
norm we denote by
$|\cdot|$. The matrix space ${\cal M}$ we assume to be
a nuclear space
densely topologically embedded in $\cal T$.
Furthermore, we assume that
${\cal M}$ is the projective limit of a family of
Hilbert spaces
$({\cal T}_p)_{p \in {\cal I}}$ parameterized by
elements of an
arbitrary set of indices ${\cal I}$ such that ${\cal
M} = \cap_{p \in
{\cal I}} {\cal T}_p$ is dense in all ${\cal T}_p$ and
that for all
$p_1, p_2 \in {\cal I}$ there
exists $p \in {\cal I}$ such that ${\cal T}_p \subset
{\cal
T}_{p_1}$ and ${\cal T}_p \subset {\cal T}_{p_2}$ and
the embeddings
are of Hilbert-Schmidt class. Let ${\cal T}_{-p}$ be
the
dual of ${\cal T}_{p}$ w.r.t.~the Hilbert
space ${\cal T}$ (if it makes no sense to consider
$-p$ of $p \in
{\cal I}$ this notation is symbolically and means that
${\cal T}_{-p}$
is the negative space corresponding to ${\cal T}_{p}$,
see
\cite{BeKo88}). Then ${\cal M}^\prime$ is the
inductive limit of the
spaces $({\cal T}_{-p})_{p \in {\cal I}}$. We denote by
$\langle \cdot,\cdot \rangle$ the dual pairings
between ${\cal T}_p$ and
${\cal T}_{-p}$ and between ${\cal M}$ and ${\cal
M}^\prime$
given by the extension of the inner product
$(\cdot,\cdot)$ on $\cal T$. We preserve this notation
for tensor
powers of these spaces.

Additionally, we introduce the notion of symmetric
tensor power of a
nuclear
space. The simplest way to do this is to start from
usual symmetric tensor
powers ${\cal T}^{\hat{\otimes} n}_p, n \in {\mathbb
N}$, of Hilbert
spaces. Using the definition
\begin{eqnarray*}
{\cal M}^{\hat{\otimes} n} :=\ \stackunder{p \in {\cal
I}}{\rm prlim} {\cal
T}^{\hat{\otimes} n}_p
\end{eqnarray*}
one can prove, see e.g \cite{Pi69} and \cite{Sch71},
that ${\cal
M}^{\hat{\otimes} n}$ is a
nuclear space which is called the $n$-th symmetric
tensor power of ${\cal
M}$. The dual space ${{\cal
M}^{\prime}}^{\hat{\otimes} n}$ can be
written as
\begin{eqnarray*}
{{\cal M}^{\prime}}^{\hat{\otimes} n} =\ \stackunder{p
\in {\cal I}}{\rm
indlim} {\cal T}^{\hat{\otimes} n}_{-p}.
\end{eqnarray*}
All the results quoted above also hold for complex
spaces.

The symmetric (or Boson) Fock space $\Gamma({\cal
T})$ of ${\cal T}$ is given by the completion of
$\mathop{\oplus}\limits_{n=0}
^\infty {\cal T}_{{\mathbb C}}^{\hat{\otimes} n}$
(${\cal T}_{{\mathbb C}}^{\hat{\otimes} 0} :=
{\mathbb C} $) w.r.t.~the Hilbertian norm
\begin{eqnarray*}
\parallel{f}\parallel^2_{\Gamma({\cal T})} :=
\sum_{n=0}^\infty n! |f^{(n)}|^2, \quad {f}=
{(}f^{(n)} {)}_{n \in {\mathbb
N}_0} \in \mathop{\oplus}\limits_{n=0}
^\infty {\cal T}_{{\mathbb C}}^{\hat{\otimes} n}.
\end{eqnarray*}

\begin{example}\label{ex1}
(i) Let ${\cal T}$ be the real vector space of all
hermitian
Hilbert-Schmidt operators defined on a complex
separable Hilbert space
${\cal H}$. Choosing an orthonormal basis $(e_i)_{i
\in {\mathbb N}}$
of ${\cal H}$ we can identify a given operator $A \in
{\cal T}$ with the
matrix
\begin{eqnarray}\label{eq2}
(a_{ij})_{i, j \in \mathbb N}, \quad a_{ij} = (e_i, A
e_j)_{\cal H},
\end{eqnarray}
where $(\cdot, \cdot)_{\cal H}$ is the inner product
of ${\cal
H}$. Since the
trace of the composition of two Hilbert-Schmidt
operators is well-defined, on ${\cal T} \times {\cal T}$ we can
introduce the bilinear form
\begin{eqnarray}\label{eq1}
(A, B) \mapsto Tr(BA) \in {\mathbb R}.
\end{eqnarray}
Since elements from ${\cal T}$ are hermitian
Hilbert-Schmidt operators
the mapping (\ref{eq1}) is a real valued, positive
definite,
symmetric, bilinear form on ${\cal T} \times {\cal
T}$. Hence on
${\cal T}$ we can define the scalar product
\begin{eqnarray*}
(A, B) := Tr(BA), \quad A,B \in {\cal T}.
\end{eqnarray*}
Now by standard methods one can verify that $({\cal
T}, (\cdot,
\cdot)$) is a real separable Hilbert space.

${\cal M}$ we choose as the subspace of ${\cal T}$
such that each
operator $F \in {\cal M}$ represented as in
(\ref{eq2}) has only
finite entries. This space can be equiped with a
topology such that it
is a nuclear space with the properties required above,
see
\cite{BeKo88}. Note that if ${\cal H}$ is not finite
dimensional
then ${\cal M}$ is not a separable nuclear space,
i.e., it can not be represented as the projective
limit of a countably
set of Hilbert spaces.

\noindent
(ii) In the finite dimensional case, i.e., ${\cal H}$
is an
$N$ dimensional, $N \in {\mathbb N}$, complex Hilbert
space, we can identify
\begin{eqnarray*}
{\cal M} = {\cal T} = {\cal M}^\prime \simeq
Mat^{h}_N({\mathbb C}),
\end{eqnarray*}
where $Mat^{h}_N({\mathbb C})$ is the real Hilbert
space of Hermitian
$(N \times N)$-matrices over the field ${\mathbb C}$.
\end{example}

In order to introduce a probability measure on the
vector space ${\cal
M}^\prime$ we consider the $\sigma$-algebra ${\cal
C}_{\sigma}({\cal
M}^\prime)$ generated by cylinder sets
\begin{eqnarray*}
C^{F_1, \ldots , F_n}_{B_1, \ldots , B_n} = \bigg{\{}
A \in {\cal
M}^\prime \bigg{|} \langle A , F_1 \rangle \in B_1,
\ldots, \langle
A , F_n \rangle \in B_n \bigg{\}},
\end{eqnarray*}
where $F_i \in {\cal M}, B_i \in {\cal B}({\mathbb
R}), \,1 \le i \le
 n, \, n \in {\mathbb N}$, and  ${\cal B}({\mathbb
R})$ denotes the
 Borel $\sigma$-algebra on ${\mathbb R}$.

The standard Gaussian measure $\mu$ on $({\cal
M}^\prime,{\cal
C}_{\sigma}({\cal M}^\prime))$ is given by its
characteristic function
\begin{eqnarray*}
C(F) = \int_{{\cal M}^\prime}  \exp(i\langle A
,F \rangle )\,d{\mu(A)} =
\exp \Big{(}-{\frac{1}{2}}|F|^2 \Big{)}, \quad F \in
{\cal M},
\end{eqnarray*}
via Minlos' theorem, see e.g.~\cite{BeKo88},
\cite{Hi80} and
\cite{HKPS93}.
$({\cal M}^\prime,{\cal C}_{\sigma}({\cal M}^\prime)$,
$\mu)$ is the basic
probability space throughout this work.
The central space in our setup is the space of complex
valued functions
which are square-integrable with respect to this
measure
\begin{eqnarray*}
L^2(\mu) := L^2({\cal M}^\prime,{\cal C}_{\sigma}
({\cal
M}^\prime),\mu).
\end{eqnarray*}
Consider the space ${\cal P}({\cal M}^\prime) \subset
L^2(\mu)$ of
smooth polynomials on ${\cal M}^\prime$:
\begin{eqnarray*}
{\cal P}({\cal M}^\prime) := \Bigg{\{} \varphi
\Bigg{|}
\varphi(A) = \sum_{n=0}^N \langle A^{\otimes n},
\tilde{\varphi}^{(n)}\rangle, \tilde{\varphi}^{(n)}
\in
{\cal M}_{{\mathbb C}}^{\hat{\otimes} n}, A \in {\cal
M}^\prime,
N \in {\mathbb N}_0 \Bigg{\}}.
\end{eqnarray*}
It is well-known that any
$ \varphi \in {\cal P}({\cal M}^\prime)$ can be
written as a
smooth Wick polynomial, i.e.,
\begin{eqnarray*}
{\cal P}({\cal M}^\prime) = \Bigg{\{} \varphi \Bigg{|}
\varphi(A) = \sum_{n=0}^N \langle :A^{\otimes n}:,
\varphi^{(n)}\rangle , {\varphi}^{(n)} \in {\cal
M}_{{\mathbb C}}^
{\hat{\otimes} n}, A \in {\cal M}^\prime, N \in
{\mathbb N}_0 \Bigg{\}}
\end{eqnarray*}
and that ${\cal P}({\cal
M}^\prime)$ is dense in $L^2(\mu)$.

Smooth Wick monomials of different order are
orthogonal w.r.t.~the
standard inner product in $L^2(\mu)$.
Furthermore, we can construct Wick monomials with
kernels
${f}^{(n)} \in {\cal T}_{{\mathbb C}}^{\hat{\otimes}
n}$ in the sense of
measurable
functions by using an approximation. More precisely,
for any sequence
$({\varphi}_j^{(n)})_{j \in {\mathbb N}_{{\mathbb C}}}
\subset {\cal
M}^{\hat{\otimes} n}_{{\mathbb C}}$ which converges to
${f}^{(n)}$ in ${\cal T}_{{\mathbb C}}^{\hat{\otimes}
n}$ we have that
$(\langle :A^{\otimes n}:, \varphi_j^{(n)}\rangle)_{j
\in \mathbb N}$
is a Cauchy sequence in $L^p(\mu), p \ge 1$, see
e.g.~\cite{BeKo88}. Since the $L^p(\mu)$ spaces are
complete the
sequence $(\langle :A^{\otimes n}:,
\varphi_j^{(n)}\rangle)_{j \in
\mathbb N}$ converges and we use
$\langle :A^{\otimes n}:, f^{(n)}\rangle$,
as a formal
notation for the limiting monomial. For Wick monomials
associated to the kernels
${f}^{(n)} \in {\cal H}_{{\mathbb C}}^{\hat{\otimes}
n}$ and
${h}^{(m)} \in {\cal T}_{{\mathbb C}}^{\hat{\otimes}
m}$, $n, m \in
{{\mathbb N}}_0$, we have the following
orthogonality property:
\begin{eqnarray*}
\int_{{\cal M}^\prime} \overline{\langle :A^{\otimes
n}:,{f}^{(n)}\rangle}
\langle :A^{\otimes m}:,{h}^{(m)}\rangle
\,d{\mu (A)} =  \delta_{n,m} \, n! \,
(\overline{f^{(n)}}, h^{(n)} )
\end{eqnarray*}
($\delta_{n,m}$ is the Kronecker delta).

The well-known Segal isomorphism between $L^2(\mu)$
and $\Gamma({\cal
T})$ establishes the It\^o-Segal-Wiener chaos
decomposition of an
element $f \in L^2(\mu)$:
\begin{eqnarray*}
f & = & \sum_{n=0}^{\infty} \langle :A^{\otimes
n}:,f^{(n)}\rangle ,
\quad f^{(n)} \in {\cal T}_{{\mathbb C}}^
{\hat{\otimes} n},
\end{eqnarray*}
and its norm is given by
\begin{eqnarray*}
\parallel f \parallel_{L^2(\mu)}^2 & = &
\sum_{n=0}^{\infty}
n! \, |f^{(n)}|.
\end{eqnarray*}

\subsection{Generalized functions on matrix
spaces}\label{ss23}

For our considerations the space $L^2(\mu)$ is too
small.
A convenient way to solve this problem is to introduce
a space of test
functions in $L^2(\mu)$ and use its larger dual space.
In Gaussian
analysis there exist various triples of test and
generalized functions
with $L^2(\mu)$ as a central space, here we choose the
Hida triple
\begin{eqnarray*}
({\cal M}) \subset L^2(\mu) \subset ({\cal
M})^{\prime},
\end{eqnarray*}
where $({\cal M})$ is the space of test functions. It
is containing the space
of smooth polynomials ${\cal P}({\cal N}^{\prime})$
and is densely and
topologically embedded in $L^2(\mu)$. Its
topologically dual
$({\cal M})^{\prime}$ is the space of Hida
distributions, see
\cite{KoSa78}, \cite{Ko80b}, \cite{Ko80a}, 
\cite{KT80}, \cite{BeKo88}, \cite{PS91}, \cite{HKPS93}, and \cite{KLPSW96}.
Instead of
reproducing the construction of this triple we give a
complete and
convenient characterization via the $S$-transform.

Consider the Wick exponential
\begin{eqnarray*}
:\exp(\langle A, F\rangle ): & := & \exp
\Big{(}\langle A,
F\rangle - {\frac{1}{2}}|F|^2 \Big{)} \\
& = & \sum_{n=0}^{\infty} \frac{1}{n!}
\langle :A^{\otimes n}:, {F}^{\otimes n}\rangle ,
\quad A \in {\cal
M}^\prime, \, F \in {\cal M}.
\end{eqnarray*}
The $S$-transform of $\Phi \in ({\cal M})^{\prime}$ at
$F \in {\cal M}$
is defined as
\begin{eqnarray*}
S\Phi(F) := \langle \! \langle \Phi, :\exp(\langle
\cdot ,F\rangle):
\rangle \! \rangle,
\end{eqnarray*}
here $\langle \! \langle \cdot, \cdot \rangle \!
\rangle$ denotes the
bilinear dual pairing between $({\cal M})$ and $({\cal
M})^{\prime}$. Since $:\exp(\langle \cdot, F\rangle ):
\in ({\cal
M})$ this pairing is well-defined. Formally we can
write the
$S$-transform as an integral
\begin{eqnarray*}
S\Phi(F) = \int_{{\cal M}^{\prime}} \Phi(A)
:\exp(\langle A,
F\rangle ): \,d\mu(A), \quad \Phi \in ({\cal
M})^{\prime}.
\end{eqnarray*}
Rigorously, this representation is valid for
e.g.~square-integrable
functions. Furthermore, we denote by ${\mathbb
E}(\Phi) := \langle \!
\langle \Phi,
1 \rangle \! \rangle = S\Phi(0)$ the (generalized)
expectation of a
Hida distribution.

In order to characterize the Hida distributions we
need the definition
of a $U$-functional.

\begin{definition}
Let $U: {\cal M} \to {\mathbb C}$ be such that:

\noindent
(i) for all $F, G \in {\cal M}$ the mapping $l \mapsto
U(G +lF)$ from ${\mathbb R}$ into $\mathbb C$ has an
entire
extension to $l \in {\mathbb C}$, and

\noindent
(ii) there exist constants $K, C > 0$ and a continuous
quadratic form
$Q$ on ${\cal M}$ so that for
all $F \in {\cal M}, z \in {\mathbb C}$,
\begin{eqnarray*}
|U(zF)| \le K \exp(C|z|^2 Q(F)).
\end{eqnarray*}
Then $U$ is called a $U$-functional.
\end{definition}

This is the base of the following characterization
theorem. For
the proof we refer to \cite{PS91} and \cite{KLPSW96}.

\begin{theorem}\label{th10}
The mapping $U: {\cal M} \to {\mathbb C}$ is the
$S$-transform of an
element in $({\cal M})^{\prime}$ if and only if it is
a $U$-functional.
\end{theorem}

Theorem \ref{th10} enables us to discuss convergence
of a sequence of
generalized functions. The following
corollary is proved in \cite{HKPS93}, \cite{PS91} and
\cite{KLPSW96}.

\begin{corollary}\label{th11}
Let $(U_n)_{n \in {\mathbb N}}$ denote a sequence of
$U$-functionals such that

\noindent
(i) $(U_n(F))_{n \in {\mathbb N}}$ is a Cauchy
sequence for all $F \in {\cal
M}$, and

\noindent
(ii) there exists a continuous quadratic form $Q$ on
${\cal M}$ and $C,D
> 0$ such that:
\begin{eqnarray*}
|U_n(zF)| \le D \exp(C|z|^2 Q(F))
\end{eqnarray*}
for all $F \in {\cal M}$, $z \in {\mathbb C}$, and
almost all $n \in
{\mathbb N}$.

\noindent
Then $(S^{-1}U_n)_{n \in {\mathbb N}}$ converges
strongly in $({\cal
M})^{\prime}$.
\end{corollary}

Sometimes it is more convenient to use the so called
$T$-transform of
a Hida distribution $\Phi$. This transform is a
generalization of the
Fourier transform in finite dimensional analysis and
is defined as
\begin{eqnarray*}
T\Phi(F) := \langle \! \langle \Phi,
\exp( i \langle \cdot ,F\rangle) \rangle \! \rangle,
\quad F \in {\cal M}.
\end{eqnarray*}
$S$- and $T$-transform have entire extensions to
${\cal M}_{\mathbb C}$ and
are related by the following formula
\begin{eqnarray*}
T\Phi(F) = C(F) \, S\Phi(iF), \quad F \in {\cal
M}_{\mathbb C}.
\end{eqnarray*}
The theorems proved above are also valid for the
$T$-transform.

\section{Symplectic Feynman distributions}\label{s0815}

Recall the matrix model with the partition function
(\ref{eq003}) defined in the introduction. 
\begin{eqnarray}\label{eq1003}
Z(N, d, g) = \int_{{\mathbb R}^{2dN^2}} \exp(iS(A,B,g)) \,dA \,dB,
\quad N, d \in {\mathbb N}, \, g \in {\mathbb R},
\end{eqnarray}
where the action is
\begin{eqnarray}\label{eq1009}
S(A, B, g) = Tr(A_{\mu}B_{\mu}) + \frac{g}{2N}
Tr(A_{\mu}B_{\nu}A_{\mu}B_{\nu}), \quad A_{\mu}, B_{\mu} \!\! \in Mat^{h}_{N}({\mathbb C}). 
\end{eqnarray}
Since the integrand in formula (\ref{eq1003}) is not integrable with
respect to the Lebesgue measure on ${\mathbb R}^{2dN^2}$ we have
to give a rigorous meaning to the heuristic formula for the partition
function. 

In this section our aim is to realize the partition function for $g = 0$ as the
generalized expectation of a Hida distribution. Firstly we regularize
the integrand in (\ref{eq1003}) and afterwards we 
show that if we remove the regulation we obtain a well-defined
element in the space of Hida distributions. 

\subsection{One dimensional case}
Let us at first consider the case $N = 1$. The suitable underlying
nuclear triple is
\begin{eqnarray*}
{\cal M} = {\cal T} = {\cal M}^\prime = Mat^{h}_{1,2d} ({\mathbb C}) \simeq
{\mathbb R}^{2d}, 
\end{eqnarray*}
where $Mat^{h}_{1, 2d} ({\mathbb C}) := (Mat^{h}_{1} ({\mathbb C}))^{2d}$
is the real vector space of $2d$-vectors having in each component a
hermitian $1 \times 1$-matrices, see Example \ref{ex1}(i), (ii).
This space is isomorphic to 
${\mathbb R}^{2d}$. The corresponding Gaussian measure we denote by
$\mu_{1, d}$; using the isomorphism to ${\mathbb R}^{2d}$ the measure
$\mu_{1, d}$ is just the standard Gaussian measure on ${\mathbb R}^{2d}$.

Let us define the regularized symplectic Feynman distribution
\begin{eqnarray*}
\Phi_{1, d, \epsilon}(a,b) := (2 \pi)^d \exp \Bigg{(} - \frac{1}{2}
\Bigg{(} \binom{a}{b}, K_{\epsilon} \binom{a}{b} \Bigg{)}_{{\mathbb
R}^{2d}} \Bigg{)}, \quad \binom{a}{b} \in {\mathbb R}^{2d},
\end{eqnarray*}
where $a ,b \in {\mathbb R}^{d}$ and $(\cdot, \cdot)_{{\mathbb
R}^{2d}}$ is the standard scalar product in ${\mathbb R}^{2d}$. The
operator $K_{\epsilon}$ is defined as
\begin{eqnarray*}
K_{\epsilon} = - \Bigg{(} (1 - \epsilon) Id_{{\mathbb R}^{2d}} 
+i \Bigg{(} \begin{array}{cc} 0 & Id_{{\mathbb R}^{d}}  \\
Id_{{\mathbb R}^{d}}  & 0 \end{array} \Bigg{)} 
\Bigg{)}, \quad \epsilon > 0.
\end{eqnarray*}
Obviously, $\Phi_{1, d, \epsilon} \in L^1(\mu_{1, d})$ and its
expectation is given by
\begin{eqnarray*}
{\mathbb E}_{\mu_{1, d}}(\Phi_{1, d, \epsilon}) 
= \int_{{\mathbb R}^{2d}} \exp \Bigg{(} - \frac{\epsilon}{2} \Bigg{(}
a_{\mu} a_{\mu} + b_{\mu} b_{\mu} \Bigg{)} + i a_{\mu} b_{\mu}
\Bigg{)} \,da \,db.  
\end{eqnarray*}
We are interested in the limit $\epsilon \to 0$ of $\Phi_{1, d,
\epsilon}$. There is no hope to find this limit in $L^1(\mu_{1,
d})$. Since for $\epsilon$ small enough the regularized
symplectic Feynman distribution $\Phi_{1, d, \epsilon}$ are not elements
from $L^2(\mu_{1,d})$ it is even a priori not clear if they are Hida
distributions. 

\begin{lemma}\label{le7001}
The regularized symplectic Feynman distributions $\Phi_{1, d,
\epsilon}, \epsilon > 0$, are Hida distributions.
\end{lemma}  
{\bf Proof:}   
Since $\Phi_{1, d, \epsilon} \in L^1(\mu_{1, d})$ we can take its $T$-transform 
at $(f, g) \in {\mathbb R}^{2d}$:
\begin{eqnarray*}
T\Phi_{1, d, \epsilon}(a,b) 
& = & (2 \pi)^d \int_{{\mathbb R}^{2d}} \exp (i(a_{\mu} f_{\mu} + b_{\mu}
g_{\mu})) \Phi_{1, d, \epsilon}(a,b) 
\,d\mu_{1, d}(a,b) \\
& = & \Bigg{(}\frac{2 \pi}{\sqrt{\epsilon^2 + 1}} \Bigg{)}^d 
\exp \Bigg{(} - \frac{1}{2}
\Bigg{(} \binom{f}{g}, \Bigg{(}K_{\epsilon} + Id_{{\mathbb R}^{2d}} \Bigg{)}^{-1}
\binom{f}{g} \Bigg{)}_{{\mathbb R}^{2d}} \Bigg{)} \\
& = & \Bigg{(}\frac{2 \pi}{\sqrt{\epsilon^2 + 1}} \Bigg{)}^d 
\exp \Bigg{(} -\frac{\epsilon (f_{\mu} f_{\mu} + g_{\mu} g_{\mu}) +
2if_{\mu}g_{\mu}}{2(\epsilon^2 + 1)} \Bigg{)}.
\end{eqnarray*}
It is easy to show that $T\Phi_{1, d, \epsilon}$ is a
$U$-functional. Thus, an application of Theorem \ref{th10} gives us that
$\Phi_{1, d, \epsilon} \in (Mat^{h}_{1, 2d} )^{\prime}$.
\hfill $\blacksquare$

Next we prove that $\Phi_{1, d, \epsilon}$ converges in $(Mat^{h}_{1,
2d} )^{\prime}$ as $\epsilon \to 0$.

\begin{proposition}\label{pr7001}
There exists a Hida distribution $\Phi_{1, d}$ such that 
\begin{eqnarray*}
\lim_{\epsilon \to 0} \Phi_{1, d, \epsilon} = \Phi_{1, d}
\end{eqnarray*}
in the strong sense. The $T$-transform of $\Phi_{1, d}$ is given by
\begin{eqnarray}\label{eq7004}
T\Phi_{1, d}(a,b) = (2 \pi)^d 
\exp (-i f_{\mu}g_{\mu}), \quad (a,b) \in {\mathbb R}^{2d}.
\end{eqnarray}
\end{proposition}  
{\bf Proof:} 
Let us choose a sequence $(\epsilon_m)_{m \in {\mathbb N}}$ which
converges to $0$. Lemma \ref{le7001} tells us that $(T\Phi_{1, d,
\epsilon_m})_{m \in {\mathbb N}}$ is a convergent sequence of
$U$-functionals. Furthermore, there exists an $M \in {\mathbb N}$ such
that
\begin{eqnarray*}
|T\Phi_{1, d, \epsilon_m} (z(f,g))| \!\! 
& \le & \!\! \Bigg{(}\frac{2 \pi}{\sqrt{\epsilon_m^2 + 1}} \Bigg{)}^d 
\exp \Bigg{(} |z|^2 \Bigg{|} \frac{\epsilon_m (f_{\mu} f_{\mu} + g_{\mu} g_{\mu}) +
2if_{\mu}g_{\mu}}{2(\epsilon_m^2 + 1)} \Bigg{|} \Bigg{)} \\
& \le & \!\! (2 \pi)^d \exp \Bigg{(} 2 |z|^2 \Bigg{|} \binom{f}{g}
\Bigg{|}_{{\mathbb R}^{2d}}^2 \Bigg{)}, \quad \forall m \ge M, \, z \in
{\mathbb C}.
\end{eqnarray*}
Therefore, the sequence $(T\Phi_{1, d,
\epsilon_m})_{m \in {\mathbb N}}$ fulfills the
properties required in Corollary \ref{th11}. Consequently,
$(T\Phi_{1, d,\epsilon_m})_{m \in {\mathbb N}}$ converges strongly in $(Mat^{h}_{1,
2d} )^{\prime}$. The limiting function we denote by
$\Phi_{1, d}$ and since strong convergence implies weak convergence
the $T$-transform of $\Phi_{1, d}$ is given by (\ref{eq7004}).  
\hfill $\blacksquare$

\subsection{Finite dimensional case}

Now we generalize these considerations to arbitrary $N \in {\mathbb
N}$. For this the appropriate basic nuclear triple is
\begin{eqnarray*}
{\cal M} = {\cal T} = {\cal M}^\prime = Mat^{h}_{N,2d} ({\mathbb C}) 
\end{eqnarray*}
and we denote by $\mu_{N, d}$ the standard Gaussian measure on this
space. For the computation of integrals w.r.t. $\mu_{N, d}$ it is
convenient to have an complete orthogonal system in $Mat^{h}_{N,d}$.
This enables the factorization of $\mu_{N, d}$ into a product of
Gaussian measures on ${\mathbb R}$. In $Mat^{h}_{N}({\mathbb C})$ equiped with
the inner product given by the trace we consider the complete orthogonal system 
\begin{eqnarray*}
COS_n = \{ e_{k,l}, 0 \le k \le l \le N\} \cup \{ iJe_{k,l}, 0 \le k < l \le N\}
\end{eqnarray*}  
where $e_{k,k}^{rs} = 1$ if $r = s = k$ and $0$ otherwise, 
and for $k \neq l$ is $e_{k,l}^{rs} = 1/2$ if $(r,s) = (k,l)$ or 
$(r,s) = (l,k)$ and $0$ otherwise. The map $J$ is given by 
multiplication with the symplectic matrix 
\begin{eqnarray*}
J = \Bigg{(} \begin{array}{cc} 0 &  Id_{{\mathbb R}^N} \\
-Id_{{\mathbb R}^N} & 0 \end{array} \Bigg{)}
\end{eqnarray*}
($Id_{\cal H}$ denotes the identity operator on the Hilbert 
space ${\cal H}$). Note that 
\begin{eqnarray*}
\Re{A}^{kl} & = & Tr(A e_{k,l}), \quad 0 \le k \le l \le N, \quad A \in
Mat^{h}_{N}({\mathbb C}), \\
\Im{A}^{kl} & = & Tr(A iJe_{k,l}), \quad 0 \le k < l \le N, 
\end{eqnarray*}
and that the system $COS_N$ is not normalized.
By standard methods the system $COS_N$ can be extended to an
orthogonal basis $COS_{N, d}$ of the product space $Mat^{h}_{N,d}
({\mathbb C})$. The $d$-vector in $COS_{N, d}$ having in the $\mu$-th
component the matrix $e_{k,l}$ or $Je_{k,l}$ and the zero matrix in the other
components we denote by $e_{\mu, k, l}$ or $Je_{\mu, k, l}$ , respectively, 
$1 \le \mu \le d$, $1 \le k \le l \le N$.  

In this setting we define the regularized symplectic Feynman
distribution 
\begin{eqnarray*}
\Phi_{N, d, \epsilon} := 2^{dN} \pi^{dN^2} \! \exp \Bigg{(} - \frac{1}{2}
\Bigg{(} \binom{A}{B}, K_{\epsilon} \binom{A}{B}
\Bigg{)}_{Mat^{h}_{N,2d}} \Bigg{)}, \, \binom{A}{B} \in  Mat^{h}_{N,2d}, 
\end{eqnarray*}
where $A, B \in  Mat^{h}_{N,d}$ and $(\cdot, \cdot)_{Mat^{h}_{N,2d}}$
is the scalar product on $Mat^{h}_{N,2d}$ given by the trace, see
Example \ref{ex1}(i), (ii) (since there is no danger of confusion we
denote from now on $Mat^{h}_{N,2d}({\mathbb C})$ by
$Mat^{h}_{N,2d}$). 
The operator $K_{\epsilon}$ is defined as
\begin{eqnarray*}
K_{\epsilon} = - \Bigg{(} (1 - \epsilon) Id_{Mat^{h}_{N,2d}} 
+i \Bigg{(} \begin{array}{cc} 0 &  Id_{Mat^{h}_{N,d}} \\
Id_{ Mat^{h}_{N,d}}  & 0 \end{array} \Bigg{)} \Bigg{)}, \quad \epsilon > 0.
\end{eqnarray*}
As in the case $N = 1$ we have $\Phi_{N, d, \epsilon} \in L^1(\mu_{N, d})$ and its
expectation is given by
\begin{eqnarray*}
{\mathbb E}_{\mu_{N, d}}(\Phi_{N, d, \epsilon}) 
= \int_{{\mathbb R}^{2dN^2}} \! \exp \Bigg{(} \! - \frac{\epsilon}{2} Tr(A_{\mu} A_{\mu} + B_{\mu}
B_{\mu}) + i Tr(A_{\mu} B_{\mu} \Bigg{)} \,dA \,dB. 
\end{eqnarray*}

\begin{lemma}\label{le2}
The regularized symplectic Feynman distributions $\Phi_{N, d,
\epsilon}, \epsilon > 0$, are Hida distributions.
\end{lemma}  
{\bf Proof:} 
We use the basis $COS_{N, d}$ of $Mat^{h}_{N, d}$ in order to evaluation the
$T$-transform of $\Phi_{N, d, \epsilon}$. Doing this, the
corresponding integral factories into a product of Gaussian measures
on ${\mathbb R}$ and we obtain 
\begin{eqnarray}\label{eq7007}
& & T\Phi_{N, d, \epsilon}(F,G) \\
& = &  2^{dN}  \Bigg{(}  
\frac{\pi}{\sqrt{\epsilon^2 + 1}} \Bigg{)}^{dN^2}
     \exp \Bigg{(}    -\frac{\epsilon Tr(F_{\mu}
F_{\mu} + G_{\mu} G_{\mu}) + 
2iTr(F_{\mu} G_{\mu})}{2(\epsilon^2 + 1)}  \Bigg{)}, \nonumber
\end{eqnarray}
$F, G \in Mat^{h}_{N,d}$. 
Obviously, $T\Phi_{N, d, \epsilon}$ is a
$U$-functional and therefore $\Phi_{N, d, \epsilon} \in (Mat^{h}_{N,
2d} )^{\prime}$ by characterization. 
\hfill $\blacksquare$

\begin{proposition}
There exists a Hida distribution $\Phi_{N, d}$ such that 
\begin{eqnarray*}
\lim_{\epsilon \to 0} \Phi_{N, d, \epsilon} = \Phi_{N, d}
\end{eqnarray*}
in the strong sense. The $T$-transform of $\Phi_{N, d}$ is given by
\begin{eqnarray*}
T\Phi_{N, d}(F,G) = 2^{dN} \pi^{dN^2} 
\exp (-i Tr(F_{\mu} G_{\mu})), \quad (F,G) \in Mat^{h}_{N,2d}.
\end{eqnarray*}
\end{proposition}  
{\bf Proof:} 
The proof of this proposition is an easy generalization of that of
Proposition \ref{pr7001}. 
\hfill $\blacksquare$

\begin{remark}\label{rm1904}
(i) Now we can realize the partition functions $Z(N, d, 0)$ as the
generalized expectation of $\Phi_{N, d}$, i.e.,
\begin{eqnarray*}
Z(N, d, 0) := {\mathbb E}_{\mu_{N, d}}(\Phi_{N, d}), \quad \Phi_{N, d}
\in (Mat^{h}_{N,2d} )^{\prime}. 
\end{eqnarray*}

\noindent
(ii) Observe, that the sequence $(\Phi_{N, d})_{N \in {\mathbb N}}$ can not
be found in one space of Hida distributions. For different $N \in
{\mathbb N}$ they are elements of different spaces of Hida
distributions for which a priori there exists no relation. Since we
are interested in the limit $N \to \infty$ next we construct a space
of Hida distributions which contains all the $\Phi_{N, d}, \, N \in
{\mathbb N}$ (up to a natural identification).
\end{remark}
   
\subsection{Infinite dimensional case}\label{ss626}

Recall the nuclear space of finite rank matrices discussed in Example
\ref{ex1}(i). We denote by ${\cal M}^{2d}$ the nuclear space which is $2d$-times the
product of this space with itself. In the same way we construct the
product space ${\cal T}^{2d}$ and obtain the nuclear triple
\begin{eqnarray}\label{eq7005}
{\cal M}^{2d} \subset {\cal T}^{2d} \subset {{\cal M}^{\prime}}^{2d}.
\end{eqnarray}
The associated standard Gaussian measure we denote by $\mu_d$. As
described in Example \ref{ex1}(i) we can represent the spaces
(\ref{eq7005}) as matrix spaces. Using this representation we define the
mapping
\begin{eqnarray*}
P_N: {{\cal M}^{\prime}}^{d} \to {\cal M}^{d}
\end{eqnarray*}  
as the projection onto the subspace of ${\cal M}^{d}$ generated by
$COS_{N,d}$ (naturally embedded into ${\cal M}^{d}$). The systems
$(COS_{N,d})_{N \in \mathbb N}$ can not only be naturally embedded into
${\cal M}^{d}$, they also can be used to constructed an orthogonal
basis of ${\cal T}^d$ which we denote by $COS_{d}$. 
We define 
\begin{eqnarray*}
\tilde{\Phi}_{N, d, \epsilon} := 2^{dN} \pi^{dN^2} \exp \Bigg{(} - \frac{1}{2}
\Bigg{(} \binom{P_NA}{P_NB}, K_{\epsilon} \binom{P_NA}{P_NB}
\Bigg{)}_{{\cal T}^{2d}} \Bigg{)}, \,\binom{A}{B} \in  {{\cal
M}^{\prime}}^{2d}, 
\end{eqnarray*}
where $A, B \in {{\cal M}^{\prime}}^{d}$ and $(\cdot, \cdot)_{{\cal T}^{2d}}$
is the scalar product on ${\cal T}^{2d}$ given by the trace, see
Example \ref{ex1}(i). The operator $K_{\epsilon}$ is defined as
\begin{eqnarray*}
K_{\epsilon} = - \Bigg{(} (1 - \epsilon) Id_{{{\cal M}^{\prime}}^{2d}} 
+i \Bigg{(} \begin{array}{cc} 0 &  Id_{{{\cal M}^{\prime}}^{d}} \\
Id_{{{\cal M}^{\prime}}^{d}} & 0 \end{array} \Bigg{)} \Bigg{)}, \quad \epsilon > 0.
\end{eqnarray*}
The regularized symplectic Feynman distributions $\tilde{\Phi}_{N, d,
\epsilon}$ and ${\Phi}_{N, d, \epsilon}$ are essentially the same. They
are only defined on different spaces and we have
\begin{eqnarray*}
{\mathbb E}_{\mu_{d}}(\tilde{\Phi}_{N, d, \epsilon}) 
& = & \int_{{\mathbb R}^{2dN^2}} \exp (- \frac{\epsilon}{2} Tr(A_{\mu} A_{\mu} + B_{\mu}
B_{\mu}) + i Tr(A_{\mu} B_{\mu}) \,dA \,dB \\
& = & {\mathbb E}_{\mu_{N,d}}(\Phi_{N, d, \epsilon}).  
\end{eqnarray*}
Hence, from now on we identify $\tilde{\Phi}_{N, d, \epsilon}$ with
${\Phi}_{N, d, \epsilon}$.  

\begin{theorem}\label{th7004}
There exists a Hida distribution $\Phi_d \in ({{\cal
M}}^{2d})^{\prime}$ such that
\begin{eqnarray*}
\lim_{N \to \infty} \lim_{\epsilon \to 0} R^{-1}{\Phi}_{N, d, \epsilon}
= \lim_{\epsilon \to 0} \lim_{N \to \infty} R^{-1} \Phi_{N, d,
\epsilon} = \Phi_d, \, R =  2^{dN}\Bigg{(}\frac{\pi}{\sqrt{\epsilon^2 +
1}} \Bigg{)}^{dN^2}, 
\end{eqnarray*}
in the strong sense.
The $T$-transform of $\Phi_{d}$ is given by
\begin{eqnarray}\label{eq11}
T\Phi_{d}(F,G) = \exp (-i Tr(F_{\mu} G_{\mu})), \quad
(F,G) \in {\cal M}^{2d}. 
\end{eqnarray}
\end{theorem}
{\bf Proof:} 
Let us choose a sequence $(\epsilon_m)_{m \in {\mathbb N}}$ which
converges to $0$. The $T$-transform of $N^{-1}{\Phi}_{N, d, \epsilon_m}$
at $(F,G) \in {{\cal M}}^{2d}$ is 
\begin{eqnarray}\label{eq6}
& & R^{-1}T\Phi_{N, d, \epsilon_m}(F,G) \nonumber \\ 
& = & \int_{{{\cal M}^{\prime}}^{2d}} \exp (i(A_{\mu} F_{\mu} +
B_{\mu} G_{\mu})) \Phi_{N, d, \epsilon_m}(A,B) \,d\mu_{d}(A,B)
\nonumber \\
& = & \exp \Bigg{(} - \frac{1}{2}
\Bigg{(} \binom{P_NF}{P_NG}, \Bigg{(}K_{\epsilon_m} + Id_{{{\cal M}^{\prime}}^{2d}} \Bigg{)}^{-1}
\binom{P_NF}{P_NG} \Bigg{)}_{{\cal T}^{2d}} \Bigg{)} \\
& \cdot & \exp \Bigg{(} - \frac{1}{2}
\Bigg{(} \binom{P^c_NF}{P^c_NG}, \binom{P^c_NF}{P^c_NG} \Bigg{)}_{{\cal
T}^{2d}} \Bigg{)}, \nonumber 
\end{eqnarray}
where $P^c_N := Id_{{{\cal M}^{\prime}}^{d}} - P_N$. The term
(\ref{eq6}), up to the normalization, coincides essentially with the
$T$-transform in the finite dimensional setting, see (\ref{eq7007}). Utilizing the
same ideas as in the proof of Proposition \ref{pr7001} we have for all
$(F,G) \in {\cal M}^{2d}$ 
\begin{eqnarray}\label{eq7008}
R^{-1}|T\Phi_{N, d, \epsilon_m}(z(F,G))| \le \exp \Bigg{(} 2 |z|^2 \Bigg{|} \binom{F}{G}
\Bigg{|}_{{\cal T}^{2d}}^2 \Bigg{)}, \quad \forall m \ge M, \, z \in {\mathbb C},
\end{eqnarray}
if we choose $M \in {\mathbb N}$ large enough. From (\ref{eq6})
together with (\ref{eq7007}) we can conclude that for all $(F,G) \in {\cal M}^{2d}$
\begin{eqnarray}\label{eq2088}
& & \lim_{N \to \infty} \lim_{m \to \infty} R^{-1}T{\Phi}_{N, d,
\epsilon}(F,G) \\
& = & \lim_{m \to \infty} \lim_{N \to \infty} R^{-1} T\Phi_{N, d,
\epsilon}(F,G) =  \exp (-i Tr(F_{\mu} G_{\mu})). \nonumber 
\end{eqnarray}
Of course, the quadratic form $Q = | \cdot |_{{\cal T}^{2d}}^2$ is
continuous on ${\cal M}^{2d}$. Thus, the estimate (\ref{eq7008}) together
with (\ref{eq2088}) and Corollary \ref{th11} proves the theorem.
                                  
\hfill $\blacksquare$

\begin{corollary}\label{co1003}
There exists a Hida distribution $\Phi_{N, d} \in ({\cal M}^{2d})^{\prime}$ such that 
\begin{eqnarray*}
\lim_{\epsilon \to 0} \Phi_{N, d, \epsilon} = \Phi_{N, d}
\end{eqnarray*}
in the strong sense. The $T$-transform of $\Phi_{N, d}$ is given by
\begin{eqnarray*}
T\Phi_{N, d}(F,G) & = & 2^{dN} \pi^{dN^2} \exp (-i Tr(P_NF_{\mu} P_NG_{\mu})) \\
& \cdot & \exp \Bigg{(} - \frac{1}{2}
\Bigg{(} \binom{P^c_NF}{P^c_NG}, \binom{P^c_NF}{P^c_NG} \Bigg{)}_{{\cal
T}^{2d}} \Bigg{)}, \quad (F,G) \in {\cal M}^{2d}. 
\end{eqnarray*}
\end{corollary}  
{\bf Proof:} 
The proof is an immediate consequence of Theorem \ref{th7004}. 
\hfill $\blacksquare$

\begin{remark}\label{rm2222}
(i) As in the finite dimensional case, see Remark \ref{rm1904}(i), we can
realize the partition 
function $Z(N, d, 0)$ as the generalized expectation of $\Phi_{N, d}$, i.e.,
\begin{eqnarray*}
Z(N, d, 0) = {\mathbb E}_{\mu_{d}}(\Phi_{N, d}), \quad \Phi_{N, d}
\in ({\cal M}^{2d})^{\prime}. 
\end{eqnarray*}

\noindent
(ii) For test functions $(F,G) \in {\cal M}^{2d}$ such that
$(P_NF, P_NG) = (F,G)$ the $T$-transform of $\Phi_{N, d}$ and $2^{dN}
\pi^{dN^2}\Phi_{d}$ coincide.
\end{remark}

\section{Proof of Theorem}\label{s777}
 
\subsection{Knots}

In order to discuss the proof of Theorem \ref{th1000} let us remind
that a knot is a smooth embedding of an oriented circle in an oriented
$3$-space ${\mathbb R}^3$. A collection of $l$ pairwise disjoint knots
is called a $l$-link. Two knots are equivalent (have the same isotopy
type) if they are equivalent under a homeomorphism of ${\mathbb R}^3$.

\begin{figure}[!htbp]
\begin{center}
\includegraphics[scale=0.75]{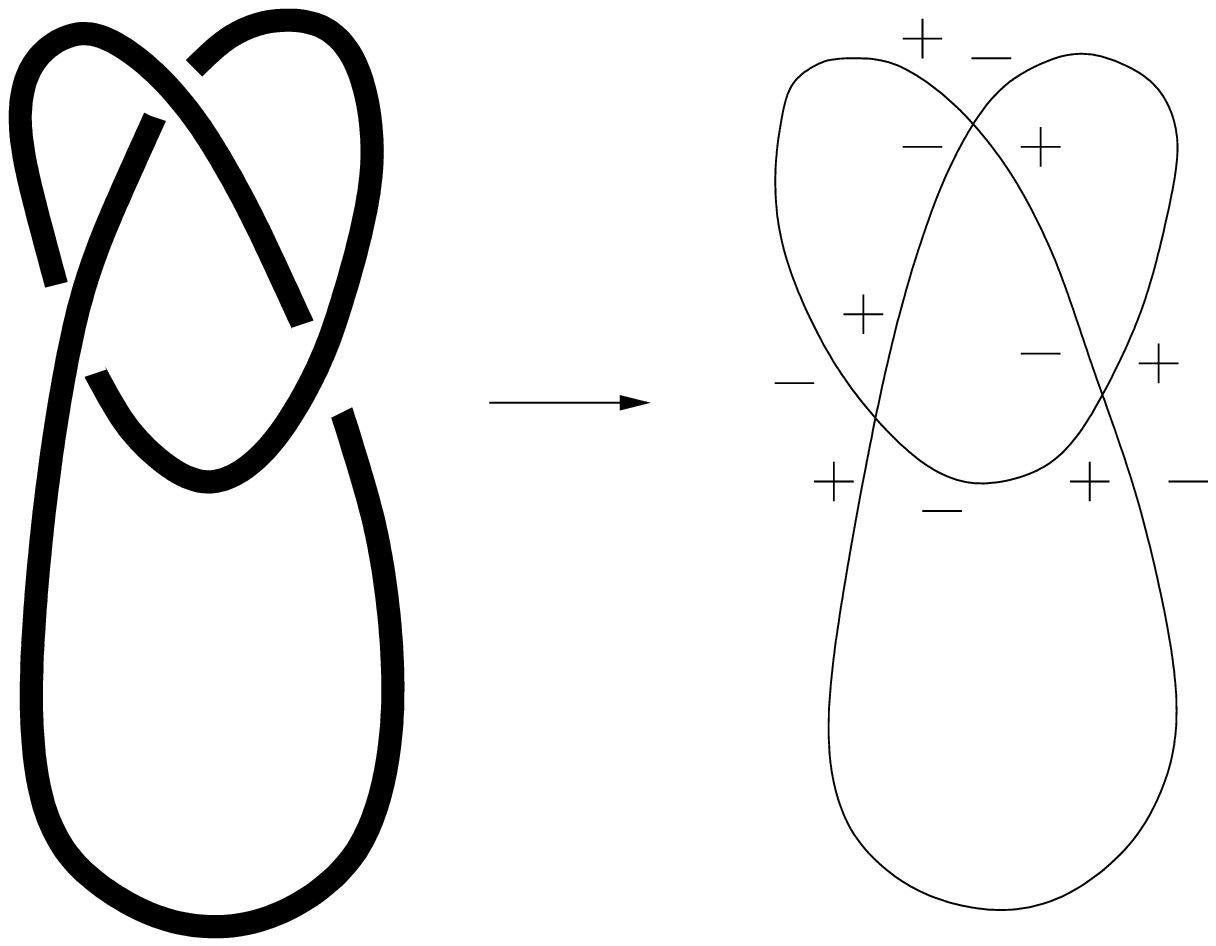}
\caption{Trefoil}
\label{fig1}
\end{center}
\begin{center}
\includegraphics{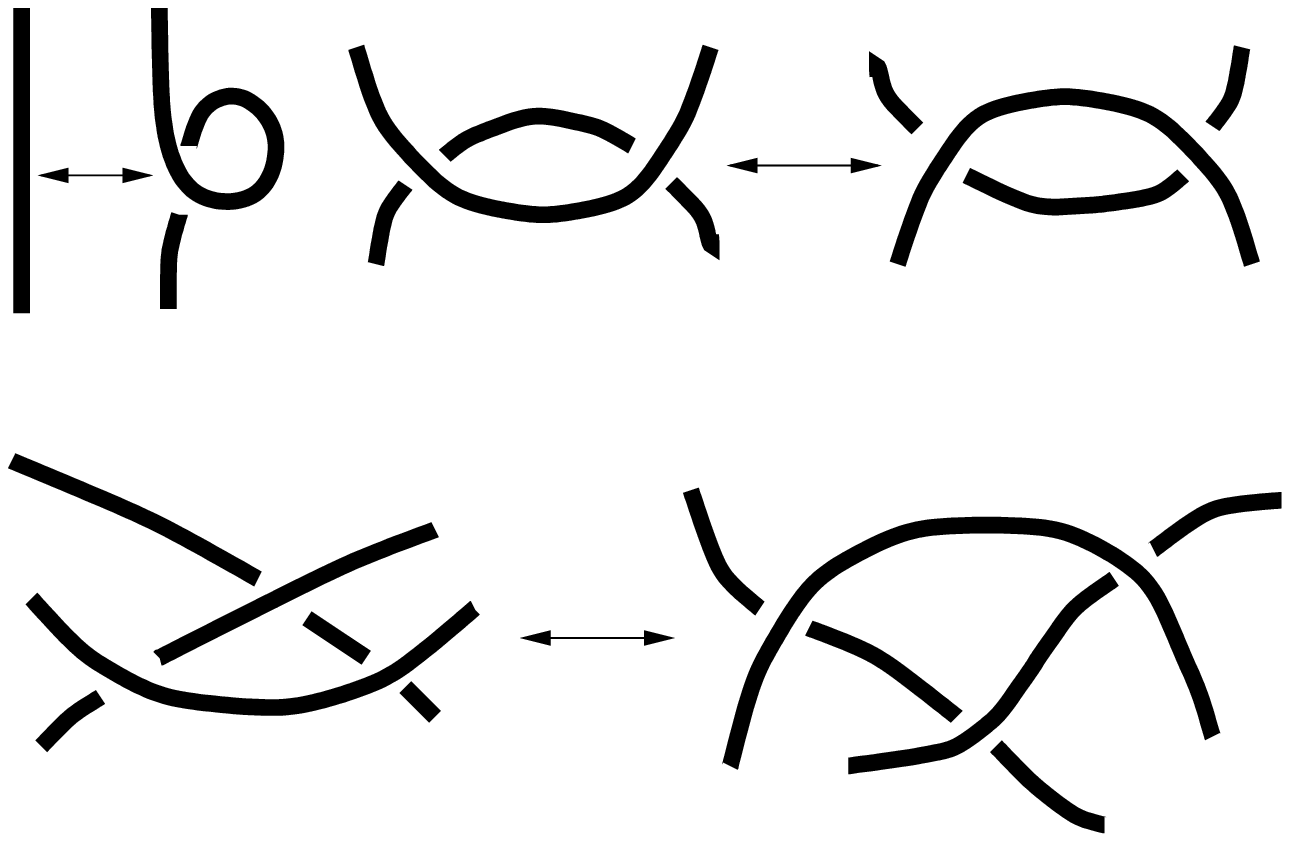}
\caption{Reidemeister moves}
\label{fig1a}
\end{center}
\end{figure}

A knot $K$ can be represented by a regular projection $\tilde{K}$ onto
the plane having at most a finite number of transverse double
points. For the plane curve $\tilde{K}$ one has to indicate which line
is up $(+)$ and which line is down $(-)$ in an intersection point, see
Figure \ref{fig1}.

In this way we get a graph on the plane which has $4$ legs in each
vertex and also has the $(+-)$ prescription. This graph is called the
knot diagram. A knot diagram is called alternating if it has
alternating $(+)$ and $(-)$ along a line. Two knot diagrams are called
Reidemeister equivalent if they define equivalent knots.
Reidemeister equivalence is generated by the three moves 
that are illustrated in Figure \ref{fig1a}.

\subsection{Wick theorem and propagators}\label{ss1860}

Now let us consider the integral (\ref{eq1003}). In order to prove
Theorem \ref{th1000}, in \cite{AV98} the authors have expanded
the partition function $Z$ into the formal
perturbation series over the coupling constant $g$ 
\begin{eqnarray}\label{eq1004}
Z =  \! \sum_{k = 0}^{\infty} \frac{1}{k!} \Big{(} \frac{ig}{N}
\Big{)}^k \! \int_{{\mathbb R}^{2dN^2}} \! (Tr(A_{\mu} B_{\nu} A_{\mu}
B_{\nu}))^k \exp(iTr(A_{\mu} B_{\mu})) \,dA \,dB. 
\end{eqnarray}
The integrals in (\ref{eq1004}) we can write as 
\begin{eqnarray}
& & \Big{(} \frac{ig}{N} \Big{)}^k \int_{{\mathbb R}^{2dN^2}}
(Tr(A_{\mu} B_{\nu} A_{\mu} 
B_{\nu}))^k \exp(iTr(A_{\mu} B_{\mu})) \,dA \,dB \label{eq4004} \\
& = & \sum_{\mu, \nu, \alpha, \beta = 1}^{d} \quad 
\sum_{q_1, i_1, j_1, l_1, m_1, n_1, o_1, p_1 = 1}^{N} \ldots
\sum_{q_k, i_k, j_k, l_k, m_k, n_k, o_k, p_k = 1}^{N} \\
& \cdot & \prod_{s = 1}^{k} \Big{(} \frac{ig}{N}
\delta_{\mu \alpha} \delta_{\nu \beta} \delta^{i_s j_s}
\delta^{l_s m_s} \delta^{n_s o_s} \delta^{p_s q_s}\Big{)}
\label{eq5005} \\ 
& \cdot & \int_{{\mathbb R}^{2dN^2}} \prod_{s = 1}^{k} 
\Big{(} A^{q_s i_s}_{\mu} B^{j_s l_s}_{\nu} A^{m_s
n_s}_{\alpha}B^{o_s p_s}_{\beta} \Big{)}\exp(iTr(A_{\mu} B_{\mu}))
\,dA \,dB \label{eq1005}     
\end{eqnarray}
($\delta_{i,j}$ and $\delta^{k,l}$ are Kronecker deltas). 
Our contribution to the proof of Theorem \ref{th1000} is to give a
mathematically rigorous meaning to the Gaussian integrals in (\ref{eq1005}), to
prove a Wick theorem for these integrals, and to derive explicit
formulas for the propagators, i.e., for the second moments of the
Gaussian integrals which we obtain in the Wick theorem.   
   
From now on we assume as underlying co-nuclear space ${{\cal
M}^{\prime}}^{2d}$.  
An elementary calculation gives
\begin{eqnarray*}
\Re{A}_{\mu}^{kl} & = & \Bigg{\langle} \binom{A}{B},
\binom{e_{\mu, k, l}}{0} \Bigg{\rangle}, \quad
\Re{B}_{\mu}^{kl} = \Bigg{\langle} \binom{A}{B},
\binom{0}{e_{\mu, k, l}} \Bigg{\rangle}, \quad k \le l,\\
\Im{A}_{\mu}^{kl} & = & \Bigg{\langle} \binom{A}{B},
\binom{iJe_{\mu, k, l}}{0} \Bigg{\rangle}, \quad
\Im{B}_{\mu}^{kl} = \Bigg{\langle} \binom{A}{B},
\binom{0}{iJe_{\mu, k, l}} \Bigg{\rangle}, \quad k < l,
\end{eqnarray*}
for $A, B \in {{\cal M}^{\prime}}^{d}$ and $e_{\mu, k, l}, \,
iJe_{\mu, k, l}\in COS_d$, where $1 \le \mu ,\nu  \le d$ and $k, l \in
{\mathbb N}$. Recall that $\Re{A}_{\mu}^{kl} = \Re{A}_{\mu}^{lk}$, 
$\Im{A}_{\mu}^{kl} = - \Im{A}_{\mu}^{lk}$ and $\Im{A}_{\mu}^{kk} = 0$
(here we assume $A, B \in {{\cal M}^{\prime}}^{d}$ in matrix 
representation as described in Example \ref{ex1}(i) and $COS_d$ naturally
embedded in this representation of ${{\cal M}^{\prime}}^{d}$ ).
This shows that the matrix coefficients and powers of them are elements
from ${\cal P}({{\cal M}^{\prime}}^{2d}) \subset ({\cal M}^{2d})$, see
Section \ref{ss23}. 

Now let us give a mathematically rigorous definition of the integral in
(\ref{eq1005}) suggested by Remark \ref{rm2222}(i). 

\begin{definition}\label{de1001}
Let $N, d, k \in {\mathbb N}$ and $\Phi_{N, d}$ as in Corollary
\ref{co1003}. Then we define the integral in (\ref{eq1005}) by
\begin{eqnarray}\label{eq1020}
&&\int_{{\mathbb R}^{2dN^2}} \prod_{s = 1}^{k}
A^{q_s i_s}_{\mu} B^{j_s l_s}_{\nu} A^{m_s
n_s}_{\alpha}B^{o_s p_s}_{\beta} \exp(iTr(A_{\mu} B_{\mu}))
\,dA \,dB  \nonumber \\ 
& := & \Big{\langle} \! \! \Big{\langle} \Phi_{N, d}, \prod_{s = 1}^{k} 
A^{q_s i_s}_{\mu} B^{j_s l_s}_{\nu} A^{m_s
n_s}_{\alpha}B^{o_s p_s}_{\beta}
\Big{\rangle} \! \! \Big{\rangle}. 
\end{eqnarray}
\end{definition}

\begin{remark}
Consider the formal perturbation series for the partition function $Z$
in (\ref{eq1004}). Using this series, formula
(\ref{eq4004})-(\ref{eq1005}), and Definition \ref{de1001}, we 
can give a mathematically rigorous definition of the partition
function $Z(N, d, g)$ as a formal power series in the coupling constant $g \in
{\mathbb R}$.
\end{remark}

At a first glance it seems to be a quite elaborate task to evaluate the integral
in Definition \ref{de1001}. This can be simplified be the following Wick theorem.

\begin{theorem}[Wick theorem]\label{th1003}
Let $F_1, \ldots, F_{2m}, \, m \in {\mathbb N}$, be elements from the
matrix space ${\cal M}^{2d}$ then 
\begin{eqnarray*}
& & \langle \! \langle \Phi_d, \langle \cdot, F_1 \rangle  \cdot \ldots
\cdot \langle \cdot, F_{2m} \rangle \rangle \! \rangle \\
& = &
\sum_{\mbox{pairings}} \langle \! \langle \Phi_d, \langle \cdot, F_{k_1}
\rangle \langle \cdot, F_{l_1} \rangle \rangle \! \rangle  
\cdot \ldots \cdot
\langle \! \langle \Phi_d, \langle \cdot, F_{k_m}
\rangle \langle \cdot, F_{l_m} \rangle \rangle \! \rangle,  
\end{eqnarray*} 
where $\sum_{\mbox{pairings}}$ means the sum over all $(2m)!/(2^m m!)$
ways of writing $1, \ldots,$ $2m$ as $m$ distinct (unordered) pairs
$(k_1, l_1), \ldots, (k_m, l_m)$.
\end{theorem}
{\bf Proof:} 
It is easy to show that 
\begin{eqnarray*}
-i\frac{d}{dt} \exp(it\langle \cdot, F \rangle)\Bigg{|}_{t=0} = \langle \cdot, F
\rangle, \quad F \in {\cal M}^{2d},
\end{eqnarray*} 
w.r.t.~the topology of $({\cal M}^{2d})$. 
Since $({\cal M}^{2d})$ is an algebra under multiplication and this
multiplication is continuous we can define the
point-wise product $\Phi_d \cdot  \langle \cdot, F
\rangle \in ({\cal M}^{2d})^{\prime}$
of the distribution $\Phi_d \in ({\cal M}^{2d})^{\prime}$ with the test
function $\langle \cdot, F \rangle  \in ({\cal M}^{2d})$ 
via the dual paring. Utilizing this product we find by an induction argument
\begin{eqnarray}\label{eq1008}
& & \langle \! \langle \Phi_d, \langle \cdot, F_1 \rangle  \cdot \ldots
\cdot \langle \cdot, F_{2m} \rangle \rangle \! \rangle \\
& = & (-1)^m \frac{\partial^{2m}}{\partial t_1 \ldots \partial t_{2m}}
T\Phi_d(t_1 F_1 + \ldots + t_{2m} F_{2m}) \Big{|}_{t_1 = \ldots =
t_{2m} = 0} \nonumber.
\end{eqnarray} 

Now under use of the explicit formula of $T\Phi_d$, see (\ref{eq11}), the
Wick theorem can be proved in the way as for the standard Gaussian
measure (the essential feature of the $T\Phi_d$ is that it is an
exponential of a quadratic form).

\hfill $\blacksquare$

\begin{remark}\label{rm4711}
Since the indices of the matrices $A$ and $B$ in Definition
\ref{de1001} are less or equal to $N$ we can identify 
\begin{eqnarray*}
\Big{\langle} \! \! \Big{\langle} \! \Phi_{N, d}, \! \prod_{s = 1}^{k} 
\! A^{q_s i_s}_{\mu} B^{j_s l_s}_{\nu} A^{m_s
n_s}_{\alpha}B^{o_s p_s}_{\beta}
\!\Big{\rangle} \! \! \Big{\rangle} \! =
\! 2^{dN} \! \pi^{dN^2} \! \Big{\langle} \! \! \Big{\langle} \!
\Phi_{d}, \! \prod_{s = 1}^{k}  
\! A^{q_s i_s}_{\mu} B^{j_s l_s}_{\nu} A^{m_s
n_s}_{\alpha}B^{o_s p_s}_{\beta}
\!\Big{\rangle} \! \! \Big{\rangle}\!. 
\end{eqnarray*}
This is an easy consequence of formula (\ref{eq1008}) and Remark \ref{rm2222}(ii).
\end{remark}

Now the evaluation of the integrals in Definition \ref{de1001} reduces
to determination of the propagators
\begin{eqnarray}\label{eq1006}
<\!A_{\mu}^{kl} B_{\nu}^{mn}\!> := \langle \! \langle \Phi_d, 
A_{\mu}^{kl} B_{\nu}^{mn} \rangle \! \rangle. \label{eq1007}
\end{eqnarray}
On the formal level can interpret 
\begin{eqnarray}\label{eq1016}
<\!A_{\mu}^{kl} B_{\nu}^{mn}\!> 
= \frac{\int A_{\mu}^{kl} B_{\nu}^{mn}
\exp(iTr(A_{\mu}B_{\mu})) \,dA \,dB} 
{\int \exp(iTr(A_{\mu}B_{\mu})) \,dA \,dB}. 
\end{eqnarray}
Analogously, we can define $<\!A_{\mu}^{kl} A_{\nu}^{mn}\!>$ and
$<\!B_{\mu}^{kl} B_{\nu}^{mn}\!>$.  

\begin{theorem}\label{pr1000}
For the propagators we have the following identities:
\begin{eqnarray*}
<\!A_{\mu}^{kl} B_{\nu}^{mn}\!> \, = \, i \delta_{\mu \nu} \delta^{kn}
\delta^{lm}, \quad <\!A_{\mu}^{kl} A_{\nu}^{mn}\!> \, = \, <\!B_{\mu}^{kl}
B_{\nu}^{mn}\!> \, = \, 0.
\end{eqnarray*}
\end{theorem}
{\bf Proof:} 
Let $k < l$ and $m > N$. In the proof of Theorem \ref{th1003} we
already derived the following formula:
\begin{eqnarray}\label{eq7009}
<\!A_{\mu}^{kl} B_{\nu}^{mn}\!> & = & - \frac{\partial^2}{\partial t_1
\partial t_2} T\Phi_d \Bigg{(} t_1 \binom{e_{\mu, k, l}}{0} 
+ t_2 \binom{0}{e_{\nu, n, m}} \Bigg{)} \Bigg{|}_{t_1 = t_2 = 0}
\nonumber \\ 
& + & \frac{\partial^2}{\partial t_1
\partial t_2} T\Phi_d \Bigg{(} t_1 \binom{iJe_{\mu, k,l}}{0} 
- t_2 \binom{0}{iJe_{\nu, n, m}} \Bigg{)} \Bigg{|}_{t_1 = t_2 = 0}
\nonumber \\
& - & i\frac{\partial^2}{\partial t_1
\partial t_2} T\Phi_d \Bigg{(} t_1 \binom{e_{\mu, k, l}}{0} 
- t_2 \binom{0}{iJe_{\nu, n, m}} \Bigg{)} \Bigg{|}_{t_1 = t_2 = 0}
\nonumber \\ 
& - & i\frac{\partial^2}{\partial t_1
\partial t_2} T\Phi_d \Bigg{(} t_1 \binom{iJe_{\mu, k,l}}{0} 
+ t_2 \binom{0}{e_{\nu, n, m}} \Bigg{)} \Bigg{|}_{t_1 = t_2 = 0} \\ 
& = & - \frac{\partial^2}{\partial t_1 \partial t_2} 
\exp \Bigg{(} - \frac{i}{2} t_1 t_2 \delta_{\mu \nu} \delta^{kn} \delta^{lm}
\Bigg{)} \Bigg{|}_{t_1 = t_2 = 0} \nonumber \\
& + & \frac{\partial^2}{\partial t_1 \partial t_2} 
\exp \Bigg{(} + \frac{i}{2} t_1 t_2 \delta_{\mu \nu} \delta^{kn} \delta^{lm}
\Bigg{)} \Bigg{|}_{t_1 = t_2 = 0} - 0 - 0 \nonumber \\
& = & i \delta_{\mu, \nu} \delta^{kn} \delta^{lm} \nonumber.
\end{eqnarray}
For $k \ge l$ and $m \le N$ we obtain the same. Having the formula
(\ref{eq7009}) in hand on easily proves that in all other cases the
propagators are equal to zero.   
\hfill $\blacksquare$

The computation of a contribution to the partition function $Z$ from
the $k$-th order can now be done by summing up Kronecker delta functions, see
(\ref{eq4004})-(\ref{eq1005}). A systematic formalism to do this is 
the Feynman diagram technique.

\subsection{Feynman diagrams}

\begin{figure}[!htbp]
\setlength{\unitlength}{0.00087489in}
\begingroup\makeatletter\ifx\SetFigFont\undefined%
\gdef\SetFigFont#1#2#3#4#5{%
  \reset@font\fontsize{#1}{#2pt}%
  \fontfamily{#3}\fontseries{#4}\fontshape{#5}%
  \selectfont}%
\fi\endgroup%
{\renewcommand{\dashlinestretch}{30}
\begin{picture}(5154,2715)(0,-10)
\path(1485,2205)(1485,495)
\path(630,1350)(2340,1350)
\path(1305,495)(1305,1170)(630,1170)
\path(2340,1170)(1665,1170)(1665,495)
\path(1665,2205)(1665,1530)(2340,1530)
\put(2500,1305){\makebox(0,0)[lb]{\smash{{{\SetFigFont{12}{14.4}{\rmdefault}{\mddefault}{\updefault}$+$}}}}}
\put(2835,1305){\makebox(0,0)[lb]{\smash{{{\SetFigFont{12}{14.4}{\rmdefault}{\mddefault}{\updefault}$\alpha$}}}}}
\put(3240,1305){\makebox(0,0)[lb]{\smash{{{\SetFigFont{12}{14.4}{\rmdefault}{\mddefault}{\updefault}$\sim$}}}}}
\put(3670,1305){\makebox(0,0)[lb]{\smash{{{\SetFigFont{12}{14.4}{\rmdefault}{\mddefault}{\updefault}$\frac{ig}{N}
\delta_{\mu \alpha} \delta_{\nu \beta} \delta^{i_s j_s}
\delta^{l_s m_s} \delta^{n_s o_s} \delta^{p_s q_s}$}}}}}
\put(360,1305){\makebox(0,0)[lb]{\smash{{{\SetFigFont{12}{14.4}{\rmdefault}{\mddefault}{\updefault}$+$}}}}}
\put(100,1305){\makebox(0,0)[lb]{\smash{{{\SetFigFont{12}{14.4}{\rmdefault}{\mddefault}{\updefault}$\mu$}}}}}
\path(630,1530)(1305,1530)(1305,2205)
\put(495,1575){\makebox(0,0)[lb]{\smash{{{\SetFigFont{12}{14.4}{\rmdefault}{\mddefault}{\updefault}$i_s$}}}}}
\put(1450,0){\makebox(0,0)[lb]{\smash{{{\SetFigFont{12}{14.4}{\rmdefault}{\mddefault}{\updefault}$\beta$}}}}}
\put(495,1050){\makebox(0,0)[lb]{\smash{{{\SetFigFont{12}{14.4}{\rmdefault}{\mddefault}{\updefault}$q_s$}}}}}
\put(1135,2250){\makebox(0,0)[lb]{\smash{{{\SetFigFont{12}{14.4}{\rmdefault}{\mddefault}{\updefault}$j_s$}}}}}
\put(1710,2250){\makebox(0,0)[lb]{\smash{{{\SetFigFont{12}{14.4}{\rmdefault}{\mddefault}{\updefault}$l_s$}}}}}
\put(1450,2580){\makebox(0,0)[lb]{\smash{{{\SetFigFont{12}{14.4}{\rmdefault}{\mddefault}{\updefault}$\nu$}}}}}
\put(1430,2330){\makebox(0,0)[lb]{\smash{{{\SetFigFont{12}{14.4}{\rmdefault}{\mddefault}{\updefault}$-$}}}}}
\put(2385,1575){\makebox(0,0)[lb]{\smash{{{\SetFigFont{12}{14.4}{\rmdefault}{\mddefault}{\updefault}$m_s$}}}}}
\put(2385,1050){\makebox(0,0)[lb]{\smash{{{\SetFigFont{12}{14.4}{\rmdefault}{\mddefault}{\updefault}$n_s$}}}}}
\put(1710,390){\makebox(0,0)[lb]{\smash{{{\SetFigFont{12}{14.4}{\rmdefault}{\mddefault}{\updefault}$o_s$}}}}}
\put(1135,390){\makebox(0,0)[lb]{\smash{{{\SetFigFont{12}{14.4}{\rmdefault}{\mddefault}{\updefault}$p_s$}}}}}
\put(1430,280){\makebox(0,0)[lb]{\smash{{{\SetFigFont{12}{14.4}{\rmdefault}{\mddefault}{\updefault}$-$}}}}}
\end{picture}
}
\caption{Feynman diagram corresponding to vertex functions}
\label{fig4}
\end{figure}
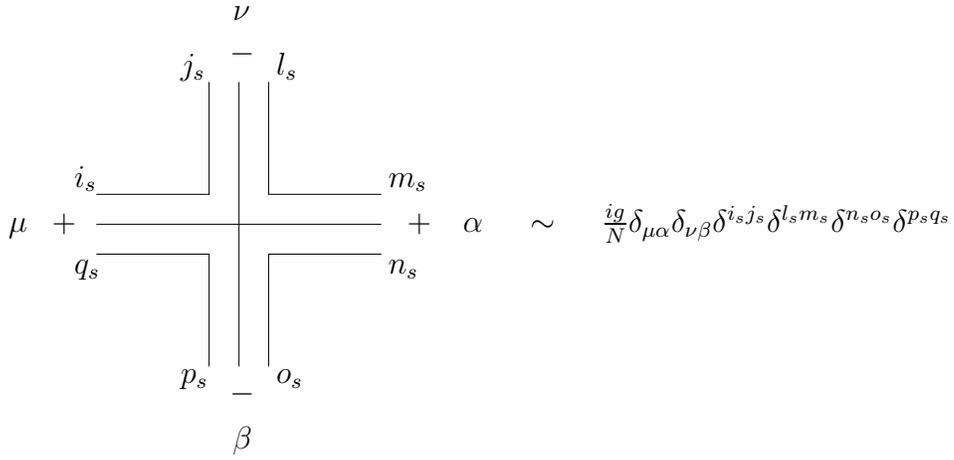

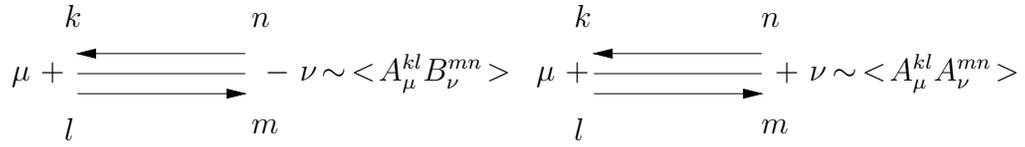
\begin{figure}[!htbp]
\setlength{\unitlength}{0.00077489in}
\begingroup\makeatletter\ifx\SetFigFont\undefined%
\gdef\SetFigFont#1#2#3#4#5{%
  \reset@font\fontsize{#1}{#2pt}%
  \fontfamily{#3}\fontseries{#4}\fontshape{#5}%
  \selectfont}%
\fi\endgroup%
{\renewcommand{\dashlinestretch}{30}
\begin{picture}(6377,969)(0,-10)
\path(495,504)(1620,504)
\path(495,369)(1620,369)
\blacken\path(1500.000,339.000)(1620.000,369.000)(1500.000,399.000)(1500.000,339.000)
\blacken\path(4090.000,669.000)(3970.000,639.000)(4090.000,609.000)(4090.000,669.000)
\path(3970,639)(5095,639)
\path(3970,504)(5095,504)
\path(3970,369)(5095,369)
\blacken\path(4975.000,339.000)(5095.000,369.000)(4975.000,399.000)(4975.000,339.000)
\put(405,819){\makebox(0,0)[lb]{\smash{{{\SetFigFont{12}{14.4}{\rmdefault}{\mddefault}{\updefault}$k$}}}}}
\put(405,54){\makebox(0,0)[lb]{\smash{{{\SetFigFont{12}{14.4}{\rmdefault}{\mddefault}{\updefault}$l$}}}}}
\put(1665,819){\makebox(0,0)[lb]{\smash{{{\SetFigFont{12}{14.4}{\rmdefault}{\mddefault}{\updefault}$n$}}}}}
\put(1665,99){\makebox(0,0)[lb]{\smash{{{\SetFigFont{12}{14.4}{\rmdefault}{\mddefault}{\updefault}$m$}}}}}
\blacken\path(615.000,669.000)(495.000,639.000)(615.000,609.000)(615.000,669.000)
\path(495,639)(1620,639)
\put(50,459){\makebox(0,0)[lb]{\smash{{{\SetFigFont{12}{14.4}{\rmdefault}{\mddefault}{\updefault}$\mu$ $+$}}}}}
\put(5185,459){\makebox(0,0)[lb]{\smash{{{\SetFigFont{12}{14.4}{\rmdefault}{\mddefault}{\updefault}$+$ $\nu$}}}}}
\put(1755,459){\makebox(0,0)[lb]{\smash{{{\SetFigFont{12}{14.4}{\rmdefault}{\mddefault}{\updefault}$-$ $\nu$}}}}}
\put(2135,459){\makebox(0,0)[lb]{\smash{{{\SetFigFont{12}{14.4}{\rmdefault}{\mddefault}{\updefault}$\sim \,<\!A_\mu^{kl}B_\nu^{mn}\!>$}}}}}
\put(5570,459){\makebox(0,0)[lb]{\smash{{{\SetFigFont{12}{14.4}{\rmdefault}{\mddefault}{\updefault}$\sim \,<\!A_\mu^{kl}A_\nu^{mn}\!>$}}}}}
\put(3835,819){\makebox(0,0)[lb]{\smash{{{\SetFigFont{12}{14.4}{\rmdefault}{\mddefault}{\updefault}$k$}}}}}
\put(3835,54){\makebox(0,0)[lb]{\smash{{{\SetFigFont{12}{14.4}{\rmdefault}{\mddefault}{\updefault}$l$}}}}}
\put(5095,819){\makebox(0,0)[lb]{\smash{{{\SetFigFont{12}{14.4}{\rmdefault}{\mddefault}{\updefault}$n$}}}}}
\put(5095,99){\makebox(0,0)[lb]{\smash{{{\SetFigFont{12}{14.4}{\rmdefault}{\mddefault}{\updefault}$m$}}}}}
\put(3580,459){\makebox(0,0)[lb]{\smash{{{\SetFigFont{12}{14.4}{\rmdefault}{\mddefault}{\updefault}$\mu$ $+$}}}}}
\end{picture}
}
\caption{Feynman diagrams corresponding to propagators}
\label{fig3}
\end{figure}

Let us interpret (\ref{eq4004})-(\ref{eq1005}) in terms of the Feynman
diagram technique. Each factor in the product (\ref{eq5005}) is a
vertex function. These vertex functions are identified with vertices
as illustrated in Figure \ref{fig4}.
The propagators $<\!A_{\mu}^{kl} B_{\nu}^{mn}\!>$ and  $<\!A_{\mu}^{kl}
A_{\nu}^{mn}\!>$ can be represented by triple lines, see Figure
\ref{fig3}.   
Each line corresponds to the separate propagation of its two
indices. The middle line carries a Greek index $\mu, \nu, \ldots$ and
the external lines carry Latin indices. The matrix $A_{\mu}$ corresponds to
$+$ and $B_{\mu}$ corresponds to $-$. 

\begin{figure}[!htbp]
\begin{center}
\includegraphics{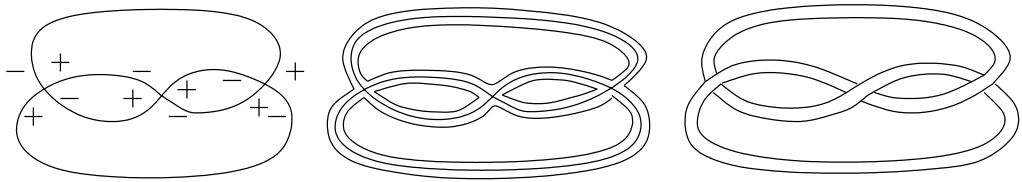}
\caption{Third order graph (trefoil) in triple line representation}
\label{fig5}
\end{center}
\end{figure}

To compute a contribution to the partition function $Z$ from the
$k$-th order of perturbation theory we draw diagrams with the 
$k$ vertices from the $k$-th order of perturbation theory, see
(\ref{eq5005}). Then we have to connect the endpoints of the vertices
corresponding to the propagators obtained from the Wick
theorem. Doing this we obtain $(4k)!/(2^{2k} (2k)!)$ diagrams, if in
such a diagram all propagators are different from zero then we call
the diagram connected. 
In Figure \ref{fig5} this is illustrated for $k = 3$ and,
additionally, it is explained how to interpret a planar connected Feynman
diagram as knot diagram.
  
The known connection between planarity and the large $N$ limit
\cite{Ho74} is based on the Euler theorem. A general Feynman diagram
consists of $P$ propagators, $V$ vertices, and $C$ closed loops of
Latin indices. The contribution of a connected diagram to the partition
function is proportional to $(g/N)^V N^C = g^V N^{C-V}$. This is clear
upon analyzing the explicit formula for the perturbation series
of the partition function $Z$, see (\ref{eq1004}) and
(\ref{eq4004})-(\ref{eq1005}) ($C$ closed loops of 
Latin indices in a connected Feynman diagram cause, in summing up the
Kronecker functions in  
(\ref{eq4004})-(\ref{eq1005}), a contribution $N^C$). Note that a 
contribution to the $k$-th order of perturbation theory may consist of
one up to $k$ connected diagrams. For a connected 
diagram one has $P = 2V$. Each closed
loop of Latin index may be considered as a face of a polyhedron and
the Euler relation reads $V - P + C = 2 - 2p$ where $p = 0, 1, \ldots$
is the number of holes of the surface on which the polyhedron is drawn
(genus of the Riemannian surface). Therefore $C - V = 2 -2p$ and the
contribution of the diagram is proportional to $g^V N^{2 - 2p}$. We
obtain that the principal contribution comes from the planar diagrams
with $p = 0$. Now in a similar vine one gets that if a planar diagram
has $l$ closed loops with Greek indices (i.e., one has a
$l$-link) then the contribution is of the diagram is proportional to
$d^l$. In particular, the contribution of the knot diagrams $l = 1$
is proportional to $d$. Hence, we obtain the following formal
expansion
\begin{eqnarray}\label{eq357}
\ln Z(N, d, g) = N \ln 2 + d N^2 \ln \pi + N^2 \sum_{l = 1}^{\infty}
\sum_{p = 0}^{\infty} F_{l, p}(g) N^{-2p} d^l.
\end{eqnarray}
Note that $Z(N, d, 0)$ is real and therefore $\ln Z(N, d, g)$ is
uniquely defined as the formal series over $g$. The summand $N \ln 2 +
d N^2 \ln \pi$ has its origin in the relation $\Phi_{N, d}= 2^{dN}
\pi^{dN^2}\Phi_{d}$, see Remark \ref{rm2222}(ii) and Remark \ref{rm4711}, 
which reflects the fact that the partition function $Z$ is not normalized. This
normalizing constant has not taken into account in \cite{AV98}.
Having derived (\ref{eq357}) except for this summand in \cite{AV98} the
authors have concluded the statement of Theorem \ref{th1000}. As
generating function for the alternating knot diagrams they have found
$F(g) = F_{1, 0}(g)$. Using (\ref{eq357}) we obtain $F(g) = \ln \pi +
F_{1, 0}(g)$. Consider $F_{1, 0}$ in its formal series
representation in the coupling constant $g$. Then the source of the
coefficient of the $k$-th power in $g$ are the planar
Feynman diagrams with $k$ vertices and one closed loop of Greek
indices. Hence, $F_{1, 0}$ can be interpreted as the generating function for
the alternating knot diagrams.

\begin{remark}
(i) One has a similar proposition for all (not only for alternating) 
knot diagrams if one takes the following action:
\begin{eqnarray}\label{eq5060}
S(A, B, g)
& = & Tr(A_{\mu}A_{\mu}) + Tr(B_{\mu}B_{\mu}) \nonumber \\ 
& + & Tr(A_{\mu}B_{\mu}) + \frac{g}{2N}
Tr(A_{\mu}B_{\nu}A_{\mu}B_{\nu}), 
\end{eqnarray}
$A, B \in Mat^{h}_{N}({\mathbb C})$. 

\begin{figure}[!htbp]
\begin{center}
\includegraphics{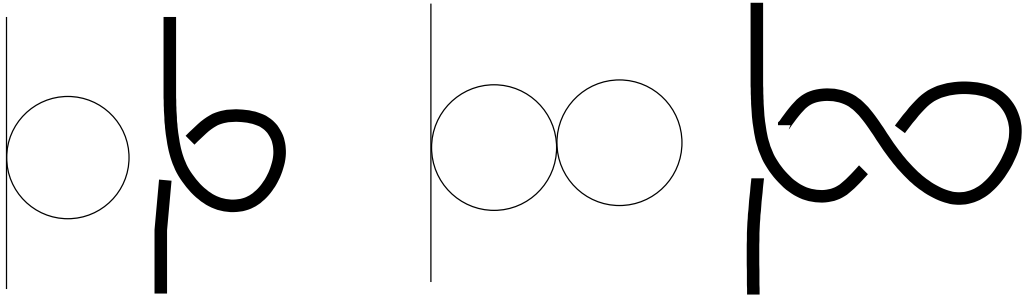}
\caption{Tadpole diagrams}
\label{fig10}
\end{center}
\end{figure}
  
\noindent
(ii) The first Reidemeister move in Figure \ref{fig1a} admits a natural 
interpretation in the Feynman diagram technique for the action 
(\ref{eq1009}). 
This move generates the so called tadpole diagrams, see
Figure \ref{fig10}, and can be 
removed by using the Wick ordered product if we use the action in the 
form    
\begin{eqnarray*}
S(A, B, g) = Tr(A_{\mu}B_{\mu}) + \frac{g}{2N}
:Tr(A_{\mu}B_{\nu}A_{\mu}B_{\nu}):_{\Phi_{N, d}}, \,\, A, B \in Mat^{h}_{N}({\mathbb C}), 
\end{eqnarray*}
where
\begin{eqnarray*}
:Tr(A_{\mu}B_{\nu}A_{\mu}B_{\nu}):_{\Phi_{N, d}}
& := & \sum_{\mu, \nu = 1}^{d}  
\sum_{j, l, m, n = 1}^{N} :A^{j l}_{\mu} B^{l, m}_{\nu} A^{m n}_{\mu}
B^{n, j}_{\nu}:_{\Phi_{N, d}} \\
& = & Tr(A_{\mu}B_{\nu}A_{\mu}B_{\nu}) + 4dn Tr(A_{\mu}B_{\mu}) + 2dN^3.
\end{eqnarray*}
Here $: \cdot :_{\Phi_{N, d}} : {\cal P}(Mat^{h}_{N,2d}) \to {\cal
P}(Mat^{h}_{N,2d})$, maps smooth polynomials to Wick ordered 
polynomials where the Wick ordering is performed w.r.t.~the
generalized function $\Phi_{N, d}$.     
\end{remark}

  


\section{Knots and M(atrix)-theory}\label{ss4711}

In \cite{AV98} the authors have drawn consequences of
Theorem \ref{th1000} concerning M-theory. Let us quote their conclusions.  

Matrix models of M-theory and superstrings are obtained by the
dimensional reduction of super Yang-Mills theory in ten dimensional
space-time to $p$ dimensions, $p = 0, 1, 2$, see \cite{BFSS96}, \cite{Pe96},
\cite{IKKT96}, and \cite{DVV97}. The bosonic part of the action in
the matrix approach to M-theory has the form
\begin{eqnarray*}
S = \int_{\mathbb R} \frac{1}{2} Tr(\dot{A}_{\mu}(t)
\dot{A}_{\mu})(t) + \frac{g}{4N} Tr([A_{\mu}A_{\nu}] [A_{\mu}A_{\nu}])
\,dt. 
\end{eqnarray*}
Here $A_{\mu}(t)$ are Hermitian $(N \times N)$-matrices over the field
$\mathbb C$ depending on time and $1 \le \mu \le d$. One has 
\begin{eqnarray}\label{eq007}
\frac{1}{2} Tr([A_{\mu}A_{\nu}] [A_{\mu}A_{\nu}]) =   
Tr(A_{\mu} A_{\nu} A_{\mu} A_{\nu}) - Tr(A_{\mu}A_{\nu}
A_{\mu}A_{\nu}).
\end{eqnarray}
The first term in (\ref{eq007}) has the form which has been discussed
in the previous sections. In the limit $d \to 0$ the principal
contribution comes from knot-like diagrams which have one loop with
Greek indices. The same reasoning one can apply to the IKKT matrix model
\cite{IKKT96} with the action
\begin{eqnarray*}
S = \frac{N}{2g} Tr([A_{\mu}A_{\nu}]^2).    
\end{eqnarray*}
If one makes the assumption on the existence of nonzero condensate \linebreak
$<\!A_{\mu}, A_{\mu}\!> \approx 1$ then one gets
\begin{eqnarray*}
S_{eff} = \frac{N}{2g} \Big{(} Tr(A_{\mu}A_{\mu}) +
Tr(A_{\mu} A_{\nu} A_{\mu} A_{\nu}) - Tr(A_{\mu}A_{\nu}
A_{\mu}A_{\nu}) \Big{)} 
\end{eqnarray*}
and one can use the described diagram technique.

Furthermore, in \cite{AV98} the authors speculate that 
Theorem \ref{th1000} and its consequences indicate that perhaps a
prototypical M-theory without matrix theory ($d = 0$) exists in the
void. The eleven dimensional M-theory could be obtained from this
prototypical M-theory by the decompactification of a point.

\subsection*{Acknowledgments}

Financial support of the DFG through the project STR
88/6-1 is gratefully acknowledged. 
Part of this work was done at CCM under the auspices of PRAXIS XXI 
(2/2.1/MAT/175/94) with support from FEDER through CITMA.
I.~Volovich is grateful to L.~Streit
for the kind hospitality at the Universidade da Madeira,
Funchal, where this work was started.

\pagebreak

\addcontentsline{toc}{section}{References}


\newcommand{\etalchar}[1]{$^{#1}$}

\end{document}